\documentclass[12pt]{amsart}
\usepackage{geometry}                
\usepackage{amsmath, amsfonts, amssymb,  latexsym, euscript, mathrsfs}

\usepackage[all]{xy}
\usepackage{hyperref}

\newcommand{\op}{{\ensuremath{\textup{op}}}}

\newcommand {\cofib} {\ensuremath{\hookrightarrow}}

\newcommand {\fibr} {\ensuremath{\twoheadrightarrow}}

\newcommand {\trivcofib} {\ensuremath{\tilde\hookrightarrow}}

\newcommand {\trivfibr} {\ensuremath{\tilde\twoheadrightarrow}}

\newcommand {\we} {\ensuremath{\tilde\rightarrow}}

\newcommand{\Lder}{\ensuremath{\mathbf{L}}} 

\newcommand{\Rder}{\ensuremath{\mathbf{R}}} 

\DeclareMathOperator{\Ho}{\textup{Ho}}

\newcommand{\cal}[1]{\ensuremath{\mathcal #1}}

\newtheorem {theorem1}{Theorem}[section]
\newtheorem {theorem}[theorem1]{Theorem}
\newtheorem {corollary}[theorem1]{Corollary}
\newtheorem {proposition}[theorem1]{Proposition}
\newtheorem {lemma}[theorem1]{Lemma}

\theoremstyle{definition}
\newtheorem {definition}[theorem1]{Definition}
\newtheorem {example}[theorem1]{Example}
\theoremstyle{remark}
\newtheorem {remark}[theorem1]{Remark}

\newcommand{\calT}{\ensuremath{\mathcal{T}} }

\newcommand{\cat}[1]{\ensuremath{\mathscr #1}}

\newcommand{\sS}{\ensuremath{\mathcal{S}} }

\newcommand{\Top}{\ensuremath{{\calT\!\text{op}\,}}}

\newcommand{\ZZ}{\ensuremath{\mathbb{Z}}}

\newcommand{\colim}{\ensuremath{\mathop{\textup{colim}}}}

\newcommand{\Sing}{\ensuremath{\mathop{\textup{Sing}}}}

\DeclareMathOperator{\Lan}{\ensuremath{\textup{Lan}}}
\DeclareMathOperator{\Ran}{\ensuremath{\textup{Ran}}}

\newcommand{\Id}{\ensuremath{\textup{Id}}}

\renewcommand{\hom}{\ensuremath{{\rm hom}}}

\newcounter{zahl}%
    {\end{list}}%

\title{A variant of a Dwyer-Kan theorem for model categories}
\author{Boris Chorny} 
\address{University of Haifa at Oranim\\
Tivon, Israel}
\email{chorny@math.haifa.ac.il}
\thanks{The first author acknowledges the support of ISF 1138/16 grant.}

\author{David White}
\address{Denison University
\\ Granville, OH, USA}
\email{david.white@denison.edu}
\thanks{The second author thanks the Center for Mathematics and Scientific Computation for supporting a visit to the University of Haifa at Oranim in 2015, when this work began.}

\date{\today}          

\usepackage[active]{srcltx}  

\begin{document}
\SelectTips{cm}{10}

\begin{abstract}
If all objects of a simplicial combinatorial model category \cat A are cofibrant, we construct the homotopy model structure on the category of small functors $\sS^{\cat A}$, where the fibrant objects are the levelwise fibrant homotopy functors, i.e., functors preserving weak equivalences. When \cat A fails to have all objects cofibrant, we construct the bifibrant-projective model structure on $\sS^{\cat A}$ and prove that it is an adequate substitute for the homotopy model structure. Next, we generalize a theorem of Dwyer and Kan, characterizing which functors $f: \cat A \to \cat B$ induce a Quillen equivalence $\sS^{\cat A} \leftrightarrows \sS^{\cat B}$ with the model structures above. We include an application to Goodwillie calculus, and we prove that the category of small linear functors from simplicial sets to simplicial sets is Quillen equivalent to the category of small linear functors from topological spaces to simplicial sets.
\end{abstract}

\maketitle

\section*{Introduction}
Homotopy functors are functors taking weak equivalences to weak equivalences. They have been a central object of interest in algebraic topology from the very beginning of the subject. W.G.~Dwyer and D.~Kan \cite{DK-diagrams} began the systematic study of the categories of homotopy functors with the theory of (what are nowadays called) relative categories (fully developed by C.~Barwick and D.~Kan, \cite{Barwick-Kan-relative}). In more detail, Dwyer and Kan ask when a map $f\colon (\cat A,U)\to (\cat B,V)$ of relative categories induces a Quillen equivalence $f^*:\sS^{\cat B,V}\to \sS^{\cat A,U}$ between the categories of homotopy functors (called restricted diagrams in \cite{DK-diagrams}) from the relative categories to the category $\sS$ of simplicial sets. Dwyer and Kan prove that $f^*$ is a Quillen equivalence if and only if the induced map of simplicial localizations $Lf\colon L(\cat A,U)\to L(\cat B,V)$  is an $r$-equivalence of simplicial categories \cite[Thm. 2.2]{DK-diagrams}.  In the current paper we formulate a version of this theorem for model categories.

Since the concept of $r$-equivalences is rarely used, especially in comparison to the concept of Dwyer-Kan equivalences, introduced in the same article, \cite{DK-diagrams},  we recall that a map $f\colon \cat A\to \cat B$ of simplicial categories is an \emph{r-equivalence} if 
\begin{itemize}
\item for every two objects $A_1, A_2\in \cat A$, the induced map 
	\[\hom_{\cat A}(A_1, A_2)\to \hom_{\cat B}(fA_1, fA_2)
	\] 
	is a weak equivalence, and 
\item every object in the ``category of components'' $\pi_0 \cat B$ is a retract of an object in the image of $\pi_0 f$.
\end{itemize}

A completely different approach to the study of the category of homotopy functors from spaces to spaces is given by Goodwillie's calculus of functors \cite{Goo:calc1, Goo:calc2, Goo:calc3}. It was noticed by W.G.~Dwyer, \cite{Dwyer-localizations}, that Goodwillie's polynomial approximation may be interpreted as a homotopical localization. This approach was reworked in terms of model categories by G.~Biedermann, the first author and O.~R\"ondigs, \cite{BCR}, in particular constructing the model category of homotopy functors on the categories of small functors from simplicial sets to simplicial sets and to spectra.

Later on, various generalizations of Goodwillie calculus to other context have appeared, \cite{unbased-calculus, Biedermann-Roendigs, Pereira}, and, hence, a natural question that arises here is the question of invariance of Goodwillie's calculus under Quillen equivalence. For example, topological spaces \Top and simplicial sets \sS are Quillen equivalent simplicial model categories. Are the model categories of small homotopy (resp., linear, $n$-excisive) functors $\sS^{\Top}$ and $\sS^\sS$ Quillen equivalent?

First we give an analog of the Dwyer-Kan theorem to model categories with all objects cofibrant: in Theorem~\ref{thm:model-for-hofunctors} we construct a model category of homotopy functors and in Theorem~\ref{Quillen-equiv} we show that a Quillen equivalence of two combinatorial model categories with all objects cofibrant gives rise to a Quillen equivalence of the categories of small functors into simplicial sets. Unfortunately this approach does not generalize further: we were not able to construct the homotopy model structure for arbitrary model categories.

The purpose of this paper is to develop a context in which the Dwyer-Kan theorem may be formulated for model categories, and then to prove that the categories of what replaces homotopy functors in our setup are equivalent if and only if the domain categories are r-equivalent. We prove this result in Theorem~\ref{main}. In particular, Example~\ref{main-example} implies together with Theorem~\ref{main} that a Quillen equivalence $\cat A \leftrightarrows \cat B$ of simplicial combinatorial model categories induces a Quillen equivalence $\sS^{\cat A}\leftrightarrows \sS^{\cat B}$ of the model categories of small homotopy functors.

The absolute version of the Dwyer-Kan theorem states that a map of simplicial categories $f\colon\cat A\to \cat B$ induces a Quillen equivalence \mbox{$\Lan_f\colon \sS^{\cat A} \leftrightarrows \sS^{\cat B} :\! f^*$} if and only if $f$ is an $r$-equivalence, \cite[Thm.~2.1]{DK-diagrams}. Luk{\'{a}}{\v{s}} Vok{\v{r}}{\'{\i}}nek generalized this result, \cite{weighted}, to categories enriched in a closed symmetric monoidal model category. The categories of homotopy functors are not discussed in his work. We give a version of the relative Dwyer-Kan theorem \cite[Thm.~2.2]{DK-diagrams} (which generalizes \cite[Thm.~2.1]{DK-diagrams}) for model categories in this paper.

As an application, we prove that the categories of small $n$-excisive functors defined on simplicial sets and on topological spaces are Quillen equivalent. More generally, given a Quillen pair such that the right adjoint preserves homotopy pushouts, we show that the model categories of $n$-excisive functors defined on this Quillen pair and taking values in simplicial sets, are Quillen equivalent. 

The paper is organized as follows. In the preliminary section we characterize which simplicial functors of simplicial combinatorial model categories induce a Quillen adjunction between the categories of small functors into simplicial sets equipped with the projective and the fibrant-projective, \cite{Duality}, model structures. In Section~\ref{bifibrant-proj} we introduce the bifibrant-projective model structure, show its existence and extend the results from the preliminary section to this new setting. Section~\ref{homotopy-model} is devoted to the study of the homotopy model structures (such that the fibrant objects are the fibrant homotopy functors) on the categories of small functors from a model category to simplicial sets. 
In order to prove the existence of the homotopy model structure, we require that the domain model category has all objects cofibrant. 
This is not a major restriction, since, as shown in \cite{Ching-Riehl,Dugger-generation}, every combinatorial model category is Quillen equivalent to one with all objects cofibrant (and in \cite{Ching-Riehl} even to one whose objects are objects of original category).
Furthermore, we show that whenever the homotopy model category exists, it is Quillen equivalent to the bifibrant-projective model structure, which exists without the requirement that all objects be cofibrant, and is a suitable replacement. Our comparison of the small functors from topological spaces to simplicial sets  with the small functors from simplicial sets to simplicial sets is  carried out in Section~\ref{motivating:example}. The homotopy model structure on $\sS^{\sS}$ was constructed in \cite{BCR} and it is Quillen equivalent to the fibrant-projective model structure. Because the Quillen model structure on $\Top$ does not have all objects cofibrant, we do not know if the homotopy model structure on $\sS^{\Top}$ exists, but the cofibrant-projective model structure on $\sS^{\Top}$ is Quillen equivalent to the fibrant-projective model structure on $\sS^{\sS}$. This means that the bifibrant-projective model structures on both categories produce Quillen equivalent model categories, as we wanted to show. 
In Section~\ref{DK} we prove our main result generalizing the Dwyer-Kan theorem. We first treat the simpler case of the homotopy model structures when all objects in the domain category are cofibrant, and then prove the general case cited above. An application to Goodwillie calculus is given in Section~\ref{Goodwillie}. We prove the Quillen equivalence of the categories of the $n$-excisive functors by localizing the Quillen equivalent categories of small functors equipped with the bifibrant-projective model structure. A tool allowing for such comparison is developed in Appendix~\ref{compare-localizations} and hopefully will be useful in other situations as well. The question of the existence of the $n$-excisive model structure is not addressed in this work, since it was considered in a number of papers, \cite{BCR, Biedermann-Roendigs, Chorny-ClassLinFun}, and the methods of localization developed there may be easily applied to the current situation.

\subsubsection*{Acknowledgements} We are grateful to Brooke Shipley and Karol Szumi{\l}o for valuable comments.

\section{Preliminaries}
In this section we recall the homotopy theory of small functors and establish some basic properties of various model categories of small functors. We assume the reader is familiar with the basics of model categories and left Bousfield localization, e.g. \cite{Hovey}, \cite{Hirschhorn}. Note that all our model categories and functors between them are simplicial, and $\hom(X,Y)$ denotes the simplicial set of morphisms from $X$ to $Y$. A model category is combinatorial if it is locally presentable and cofibrantly generated.

\begin{definition}
Let \cat A be a simplicial category. A functor $F\colon \cat A \to \sS$ is \emph{small} if it is a left Kan extension from some small subcategory. In other words, there exists a small full subcategory $i\colon \cat A' \hookrightarrow \cat A$, such that $F=\Lan_i i^*F$. We denote the category of small functors from $\cat A$ to $\sS$ by $\sS^{\cat A}$.
\end{definition}

\begin{remark}
In the book by M.~Kelly, \cite{Kelly}, small functors are called accessible, which does not correspond to modern terminology (accessible functors are functors of accessible categories preserving $\lambda$-filtered colimits for some cardinal $\lambda$), though accessible functors are always small and small functors of accessible categories are accessible. A  functor is small if and only if it is a small (weighted) colimit of representable functors, \cite[Proposition~4.83]{Kelly}. Since the category of small functors from \cat A to \sS  is cocomplete, \cite[Proposition~5.34]{Kelly}, in particular tensored over $\sS$, a colimit of the functor $G\colon (\cat A')^{\op}\to \sS^{\cat A}$  weighted by the functor $F\colon \cat A' \to \sS$ may be computed using the coend formula: $F\star_{\cat A'} G = \int^{A\in \cat A'} FA \otimes GA$, \cite [3.70]{Kelly}. 
\end{remark}

Now we would like to analyze what kind of functors are induced on the categories of small functors by an adjunction of domain categories.
 
\begin{proposition}\label{adjunction}
Let $L\colon \cat A\to \cat B$ be a simplicial accessible functor between locally presentable simplicial categories. Then there exists a pair of adjoint functors between the categories of small functors 
\[
\Lan_L\colon \sS^{\cat A} \leftrightarrows \sS^{\cat B} :\! L^*.
\] 
If in addition $L$ has a right adjoint $R$, then $\Lan_L = R^*$ is given by the precomposition with $R$.
\end{proposition}

\begin{proof}
Note that every small functor $F\in \sS^{\cat A}$ is a left Kan extension from a full small subcategory $i\colon \cat A' \hookrightarrow \cat A$. Then
\[
\Lan_L(F) = \Lan_L(\Lan_i i^*F) =  \Lan_{Li} i^*F
\]
by the transitivity property of the iterated left Kan extensions, \cite[Theorem~4.47]{Kelly}. Hence, $\Lan_L(F)\in \sS^{\cat A}$ is a small functor.

Given a representable functor $R^B = \hom_{\cat B}(B,-)$, the functor $L^*R^B = \hom(B,L-)$ is no longer representable, but it is $\lambda$-accessible if $B$ is $\lambda$-presentable and $L$ is $\lambda$-accessible. Hence it is a small functor as an accessible functor of accessible categories.

For any $G\in \sS^{\cat B}$, 
\[
L^*G = L^* \left(\int^B \hom(B,-)\otimes GB\right) = \int^B \hom(B,L-)\otimes GB
\]
is a weighted colimit of small functors, which is again small \cite[5.34]{Kelly}.

Suppose now that $L$ has a right adjoint $L\dashv R$. Then, using Yoneda's Lemma, 
\begin{align*}
\hom(R^* F, G)=\hom(R^* \left(\int^{A\in \cat A'} \hom(A,-)\otimes FA\right), G)= \\ 
\int_A\hom(\hom(A,R(-))\otimes FA, G)= \int_A\hom(FA, \hom(\hom(LA,-),G))=\\
\int_A\hom(FA, G(LA))= \int_A\hom(FA, L^*G(A))=\\
 \int_A\hom(FA, \hom(\hom(A,-),L^*G))=
\int_A\hom(\hom(A,-)\otimes FA ,L^*G)=\\
\hom \left(\int^A\hom(A,-)\otimes FA ,L^*G\right)=\hom(F, L^*G).
\end{align*}

In other words $R^* \dashv L^*$, hence $R^*=\Lan_L$.
\end{proof}

We are interested in the homotopy theory of small functors. The projective model structure (where weak equivalences and fibrations are levelwise) on the category of small functors was constructed in \cite[Theorem~3.1]{Chorny-Dwyer} for all cocomplete domain categories. The condition of cocompleteness is required to ensure that the category of small functors is complete, \cite[Corollary~3.9]{Day-Lack}.
 
\begin{proposition}\label{Quillen-map}
Given a simplicial accessible functor $f\colon \cat A \to \cat B$ of simplicial combinatorial model categories, the adjunction $\Lan_f \dashv f^*$ discussed in Proposition~\ref{adjunction} is  a Quillen pair for the projective model structures on the categories of small functors from $\cat A$ and $\cat B$ to $\sS$. 
\end{proposition}
\begin{proof}
Consider the adjunction
\[
\xymatrix{
\sS^{\cat A}
\ar@/^/[r]^{\Lan_f}
				&	\sS^{\cat B}
					\ar@/^/[l]^{f^*}
}
\]
between the two model categories of small functors equipped with the projective model structure, \cite{Chorny-Dwyer}. Let  $p\colon F\to G$ be a (trivial) fibration in $\sS^{\cat B}$. Consider the induced map $f^*p\colon f^*F \to f^*G$ in $\sS^{\cat A}$. Let $A\in \cat A$ be an arbitrary object. Then $p_{fA}\colon F(fA)\to G(fA)$ is a (trivial) fibration by assumption. Furthermore, $f^*p_A = p_{fA}$ is also a (trivial) fibration:
\[
\xymatrix{
(f^*F)(A)
\ar[r]^{f^*p_A}
\ar@{=}[d]
			&	(f^*G)(A)
				\ar@{=}[d]\\
F(fA)
\ar@{->>}[r]^{(\dir{~}\,\,\,)}_{p_{fA}}
			&	G(fA)
				\\
}
\]  

\end{proof}

The fibrant-projective model structure on the category of small functors with domain in a combinatorial model category (where weak equivalences and fibrations are levelwise in fibrant objects) was constructed in \cite[Definition~3.2]{Duality}. This is a particular case of the relative model structure, \cite[Definition~2.2]{Chorny-relative}. In the next proposition we analyze its interaction with the adjunction of Proposition \ref{adjunction}.

\begin{proposition}\label{Q-pair-fib-proj}
Given a simplicial accessible functor $f\colon \cat A \to \cat B$ of simplicial combinatorial model categories, the adjunction $\Lan_f \dashv f^*$ discussed in Proposition~\ref{adjunction} is  a Quillen pair for the fibrant-projective model structure on the categories of small functors from $\cat A, \, \cat B$ to $\sS$ if and only if $f$ preserves fibrant objects. 
\end{proposition}
\begin{proof}
The ``if'' direction follows in the same manner as Proposition \ref{Quillen-map} above.

For the ``only if'' direction, assume that $f^*$ is a right Quillen functor and we need to show that for every fibrant $A\in\cat A$ the map $p\colon fA \to \ast$ has the right lifting property with respect to any trivial cofibration $i\colon B_1 \trivcofib B_2$ in \cat B. By \cite[Prop.~9.4.3]{Hirschhorn}, it suffices to show that $(i,p)$ is a homotopy lifting-extension pair. In other words, it suffices to show that $\hom(B_2, fA)\to \hom(B_1, fA)$ is a trivial fibration of simplicial sets.

For any trivial cofibration $i\colon B_1 \trivcofib B_2$ in \cat B the induced map of representable functors $i^*\colon \hom(B_2,-)\to \hom(B_1,-)$ is a trivial fibration in the fibrant-projective model structure on $\cal S^{\cat B}$, by the SM7 axiom, \cite[Defn.~9.1.6]{Hirschhorn}. Since $f^*$ is a right Quillen functor, the map 
\[
f^*i^*\colon \hom(B_2, f-)\to \hom(B_1, f-)
\]
 is a trivial fibration in the fibrant-projective model structure on $\cal S^{\cat A}$, i.e., 
 \[
 \hom(B_2, fA)\trivfibr \hom(B_1, fA)
 \]
 is a trivial fibration of simplicial sets for all fibrant $A\in \cat A$.
\end{proof}

\section{Bifibrant-projective model structure}\label{bifibrant-proj}

Let $\cat A$ be a simplicial combinatorial model category. By analogy with the fibrant-projective, \cite[Definition~3.2]{Duality}, and the cofibrant-projective, \cite[Definition~2.1]{Blanc-Chorny}, model structures on the categories of small functors $\sS^{\cat A}$, we introduce the bifibrant-projective model structure on $\sS^{\cat A}$.

\begin{definition}\label{bifibrant-def}
Let \cat A be a simplicial combinatorial model category, and let $F, G\in \sS^{\cat A}$ be small functors. A natural transformation $f\colon F\to G$ is a \emph{bifibrant weak equivalence} (resp., \emph{bifibrant fibration}) if for all bifibrant objects $A\in \cat A$ (i.e., objects which are both fibrant and cofibrant), the induced map $f_A\colon F(A)\to G(A)$ is a weak equivalence (resp., a fibration) of simplicial sets. A natural transformation is a \emph{bifibrant cofibration} if it has the left lifting property with respect to bifibrant trivial fibrations.
\end{definition}

Next, we establish the existence of the bifibrant-projective model structure as a particular case of the relative model structure, \cite[Definition~2.1]{Chorny-relative}.

\begin{proposition}\label{bifibrant-exist}
Let \cat A be a simplicial combinatorial model category. Then the category of small functors $\sS^{\cat A}$ may be equipped with the bifibrant model structure.
\end{proposition}
\begin{proof}
We will verify the conditions of \cite[Proposition~2.8]{Chorny-relative} in order to establish the bifibrant model structure, which is also the bifibrant relative model structure in the terminology of \cite{Chorny-relative}.

The condition requiring verification is the local smallness, \cite[Definition~2.4]{Chorny-relative}, of the subcategory of bifibrant objects in the category $\cat A^{\op}$, or, dually, the solution set condition in \cat A, i.e., for every object $A\in \cat A$ we need to find a \emph{set} of bifibrant objects $\mathscr{W}_A$ such that for every bifibrant object $B$ and every map $A\to B$ there exists a (non-unique) object $W\in \mathscr{W}_A$ such that $A\to W\to B$.

For every object  $A\in \cat A$, choose a cardinal $\kappa$ large enough, so that $A$ is $\kappa$-presentable and \cat A is $\kappa$-combinatorial. Next, look at the set $\mathscr W_A'$ of $\kappa$-presentable cofibrant objects in \cat A, then put $\mathscr{W}_A=\{\widehat W \, |\, W\in \mathscr W_A' \}$, where $\widehat W$ denotes fibrant replacement.

The fat small object argument, \cite[Corollary~5.1]{fat}, shows that every cofibrant object is a $\kappa$-filtered colimit of $\kappa$-presentable cofibrant objects in the $\kappa$-combinatorial model category \cat A. It follows that every morphism $A\to B$ into a bifibrant object $B$ factors first through some $W_1\in \mathscr W_A'$. Finally, the morphism $W_1\to B$ factors through the fibrant replacement $W=\widehat{W}_1$ of $W_1$, since $B$ is fibrant: $W_1\trivcofib W \to B$.
\end{proof}

Now we need to find the conditions on a functor $f\colon \cat A \to \cat B$ between simplicial combinatorial model categories, so that the induced adjunction $Lan_f \dashv f*$ of Proposition \ref{adjunction} is a Quillen pair. 

\begin{proposition}\label{Quillen-pair-bifibrant}
Given a simplicial accessible functor $f\colon \cat A \to \cat B$ of simplicial combinatorial model categories, the adjunction $\Lan_f \dashv f^*$ discussed in Proposition~\ref{adjunction} is a Quillen pair for the bifibrant-projective model structure on the categories of small functors from $\cat A$ and $\cat B$ to $\sS$ if $f$ preserves both fibrant and cofibrant objects. 
\end{proposition}
\begin{proof}
Similar to Proposition~\ref{Quillen-map}.
\end{proof}

\begin{example}
The classical Quillen equivalence $|-|\colon \sS \leftrightarrows \Top :\! \Sing$ induces the following Quillen map of bifibrant-projective model structures
\[
{\Sing}^*\colon \sS^\sS \leftrightarrows \sS^\Top :\! |-|^*,
\]
which turn out to be the fibrant-projective and the cofibrant-projective model structures respectively. Of course, this is a very special case when the left Quillen functor preserves fibrant objects. We will have to find a way around this difficulty in order to generalize this example.
\end{example}

\section{Homotopy model structure}\label{homotopy-model}
Let \cat A be a simplicial combinatorial model category. Recall that \emph{homotopy functors} are functors preserving the weak equivalences. If it exists, the \emph{homotopy model structure} on the category of small functors $\cal S^{\cat A}$ is a localization of the projective model structure in such a way that the local objects are the projectively fibrant homotopy functors. We will only construct the homotopy model structure on the category of small presheaves $\cal S^{\cat A}$ under the additional assumption that all objects of \cat A are cofibrant.

\subsection{Localization construction} \label{precomposition}

If we localize the projective model structure on the category of small functors $\cal S^{\cat A}$ with respect to the following class of maps
\[
H_{\cat A} = \left\{\left. \hom(A_1, -) \to \hom(A_2, - ) \right| A_1\we A_2 \text{ in } \cat A  \right\},
\]
then the fibrant objects in the new model structure will be precisely the levelwise fibrant homotopy functors. The resulting model structure is the homotopy model structure on $\cal S^{\cat A}$. Since the projective model structure is not cofibrantly generated (it has a proper class of generating cofibrations, instead of a small set), and $H_{\cat A}$ is a proper class of maps, the localization techniques of Smith and Hirschhorn may not be applied. 

In the case that all objects of \cat A are cofibrant, we will use the Bousfield-Friedlander, \cite[Appendix~A]{BF}, $Q$-model structure construction further improved by Bousfield, \cite[Theorem~9.3]{Bou:telescopic} in order to obtain the left Bousfield localization of $\sS^{\cat A}$ with respect to $H_{\cat A}$. 

\begin{theorem} \label{thm:model-for-hofunctors}
Let \cat A be a simplicial combinatorial model category with all objects cofibrant. Then there exists a localization of the projective model structure on $\cal S^{\cat A}$, such that the fibrant objects are precisely the levelwise fibrant homotopy functors.
\end{theorem}

\begin{proof}
Since \cat A is a simplicial combinatorial model category, we can fix a continuous, accessible fibrant replacement functor $\mathrm{Fib_{\cat A}}\colon \cat A\to \cat A$ together with a natural transformation $\varepsilon\colon \Id_{\cat A} \to \mathrm{Fib}_{\cat A}$. These properties are required to ensure that the precomposition of $\mathrm{Fib_{\cat A}}$ with a small functor $F\colon \cat A\to \sS$ produces a small functor again.

We denote the fibrant replacement in $\cal S$ by $\widehat{(-)}$. In this case the homotopy approximation functor may be constructed very explicitly. Namely, for any small $F\colon \cat A \to \cal S$, we can put $\cal H(F) = \widehat{\mathrm{Fib}^\ast_{\cat A} F} = \widehat{F\circ \mathrm{Fib_{\cat A}}}$. It is equipped\ with the coaugmentation: $\widehat{F\varepsilon}\colon F \to \widehat{\mathrm{Fib}^\ast_{\cat A} F}$. This is a homotopy idempotent construction that takes values in homotopy functors, since weak equivalences of objects which are fibrant and cofibrant are simplicial weak equivalences, \cite[II.2.5]{Quillen}, and the latter are preserved by simplicial functors, cf. \cite[Proposition~3.3]{BCR}.  

By \cite[Proposition~4.3]{Raventos}, in any model category $\cat M$ equipped with  a homotopy idempotent functor $L\colon \cat M\to \cat M$, the class of $L$-equivalences (the maps rendered by $L$ into weak equivalences) coincides with the class of the local equivalences (the class of maps simplicially orthogonal to the $L$-local objects), therefore  $\cal H$-equivalences are precisely the local equivalences with respect to the fibrant homotopy functors. Since our construction is very simple, we can see immediately that $\cal H$-equivalences, i.e., maps rendered into projective weak equivalences by the functor $\cal H$, are precisely the fibrant-projective weak equivalences of small functors, \cite[Def.~3.2]{Duality}, i.e., the natural transformations of functors inducing weak equivalences of fibrant objects.

It remains to verify that our localization construction satisfies the conditions A1-A3 of \cite[Theorem~9.3]{Bou:telescopic}. The projective model structure on the category $\sS^{\cat A}$ of small functors is proper by \cite[Theorem~3.6]{Duality}, since \sS is a right proper model category and a strongly left proper monoidal model category, \cite[Definition~4.5]{DRO:enriched}.

A1 and A2 are satisfied by the construction of $\cal H$ and the discussion above. To verify A3 consider the pullback of a fibrant-projective weak equivalence along a projective fibration. Since \sS is right proper, the base change of  a fibrant-projective weak equivalence is a fibrant-projective weak equivalence again.

Hence the left Bousfield localization exists, and defines the $\cal H$-local model structure on the category of small functors from \cat A to $\sS$. This is the homotopy model structure, since the $\cal H$-local objects are precisely the projectively fibrant homotopy functors. In other words, $\cal H$-localization is the localization with respect to $H_{\cat A}$. 
\end{proof}


If we drop the assumption that all objects are cofibrant, we are unable to construct the left Bousfield localization of the projective model category $\sS^{\cat A}$ with local objects being precisely the homotopy functors, but we will show the existence of a homotopy idempotent (non-functorial) localization construction $Q$, such that $Q$-equivalences are precisely the $H_{\cat A}$-equivalences. 

\begin{proposition}\label{homotopy-approx}
Let \cat A be a simplicial combinatorial model category. Then for each functor $F\in \sS^{\cat A}$ there exists an $H_{\cat A}$-equivalence $F\to QF$ such that $QF$ is a homotopy functor.
\end{proposition}

\begin{proof}
Let $f\colon \cat A\to \cat A$ be a bifibrant replacement functor. Notice that the adjunction $\Lan_f \dashv f^*$ is a Quillen pair for the projective model structure on $\sS^{\cat A}$, since the right adjoint $f^*$ preserves fibrations and trivial fibrations. 

For every functor $F\in H_{\cat A}$ consider the following construction:
\[
\xymatrix{
\tilde{F}
\ar@{->>}[dd]_{\dir{~}}
\ar[rr]
\ar@{^(->}[dr]	&	& f^* \Lan_f \tilde F\\
	& Q'F
	   \ar@{->>}[ur]^{\dir{~}}
	   \ar[d]^{\dir{~}}\\
F
\ar[r]	&QF\\
}
\]
We begin with $F\in H_{\cat A}$, take its cofibrant replacement $\tilde F$ in the projective model structure and factor the unit of the adjunction $\Lan_f \dashv f^*$ into a projective cofibration followed by a projective trivial fibration. Denote the middle term of the factorization by $Q'F$ and put $QF=Q'F\coprod_{\tilde F} F$.
 
Unfortunately this construction is not functorial, since the cofibrant replacement in the projective model structure is not known to be functorial. On the other hand, it naturally extends to morphisms and can be rendered functorial on any small subcategory of $\sS^{\cat A}$.

This construction preserves weak equivalences of functors, since all the stages of the construction do. In particular, $\Lan_f$ preserves weak equivalences between cofibrant objects and $f^*$ preserves all projective weak equivalences.

Now, $QF$ is projectively weakly equivalent to $f^* \Lan_f \tilde F$, so in order to show that $QF$ is a homotopy functor it suffices to show that $f^* \Lan_f \tilde F$ is. We will show it by cellular induction, assuming that $\tilde F=\colim_{i < \lambda} F_i$, so that $F_0=\emptyset$ and $F_{i}$ is obtained from $F_{i-1}$ by attaching a cell
\[
\xymatrix{
R^{A_i}\otimes \partial\Delta^n
\ar@{^(->}[d]
\ar[r]			 & F_{i-1}
				\ar@{^(->}[d]\\
R^{A_i}\otimes \Delta^n
\ar[r]			 & F_{i}\\
}
\]
if $i$ is a successor ordinal or $F_i=\colim_{a<i}F_a$ if $i$ is a limit ordinal.

In order to compute $\Lan_f \tilde F$, notice that $\Lan_f$ commutes with colimits, and $\Lan_f R^{A_i}=R^{f(A_i)}$. In other words, $\Lan_f \tilde F$ is a  cellular complex with cells of type $R^{f(A_i)}$, i.e., represented in cofibrant objects. 

Next, we must show that $QF$ is a homotopy functor, i.e., it preserves weak equivalences. The following argument proves that a projectively equivalent functor $f^*\Lan_f \tilde F$ is a homotopy functor by cellular induction. First, note that $f^*$ preserves colimits, as a left adjoint to $\Ran_f$, which exists, in turn, by \cite{Day-Lack}, or just because the colimits in the diagrams of functors are computed levelwise. Therefore, $f^*\Lan_f \tilde F$ is a cellular construction, with cells of type $f^*R^{f(A_i)} = \hom(f(A_i), f(-))$, so it is no longer a representable functor, but is a homotopy functor. Hence, assuming for induction that $f^*\Lan_f F_{a}$ is a homotopy functor for all $a < i$, we obtain that $f^*\Lan_f F_{i}$ is also a homotopy functor, hence $f^*\Lan_f \tilde F$ is a homotopy functor as a sequential colimit of homotopy functors into $\sS$.

The last claim that we need to show is that the map $F\to QF$ is an $H_{\cat A}$-equivalence. In other words, that our construction is homotopy idempotent. We will show the equivalent statement that the map $F\to QF$ is initial in a suitable sense, i.e., we will show that in the homotopy category $\Ho(\sS^\cat A)$ the unit of the derived pair of adjoint functors $\varepsilon\colon [\tilde F]=[F]\to [QF]= \Rder f^* \Lder \Lan_f [\tilde F]$ is initial with respect to maps into homotopy functors.  

Let $H\in \sS^{\cat A}$ be a projectively fibrant homotopy functor, and let $g\colon [F]\to [H]$ be a map in the homotopy category. Notice that since $H$ preserves weak equivalences, $[H]=[f^*H]=\Rder f^* [H]$. Then by the universal property of the unit there exists a unique map $h\colon \Rder f^* \Lder \Lan_f \tilde F \to \Rder f^* [H]=[H]$ such that $g=h\varepsilon$.

Therefore, our initial construction $F\to QF$ is a homotopy localization turning every small functor into a homotopy functor, i.e., localization with respect to $S_\cat A$.
\end{proof}

\subsection{Comparison of the homotopy and the bifibrant-projective model structures}
We next prove that the homotopy model structure on $\sS^\cat A$ is Quillen equivalent to the bifibrant-projective model structure (Definition~\ref{bifibrant-def}), when both model structures exist. This means we can use the bifibrant-projective model structure as a substitute for the homotopy model structure in contexts where the homotopy model structure is not yet known to exist. We conjecture that the homotopy model structure on $\sS^\cat A$ exists even if not all objects of $\cat A$ are cofibrant. We expect that localization of class-combinatorial model categories \cite{Chorny-Rosicky-II} can be used to prove this conjecture.

For the sake of comparison, we assume in this section that the homotopy model structure exists. We note that the bifibrant-projective model structure exists whenever \cat A is combinatorial. We show now that these model structures are Quillen equivalent. 

\begin{theorem}\label{fib-proj-equiv}
Let \cat A be a simplicial combinatorial model category. Then the pair of identity functors induces a Quillen equivalence of the homotopy and the bifibrant-projective model structures.
\end{theorem}
\begin{proof}
Consider the pair of adjoint functors
\[
\Id\colon \cal S_{\text{bifib-proj}}^{\cat A} \rightleftarrows \cal S_{\text{proj}}^{\cat A} :\! \Id,
\]
where the left adjoint is pointing from left to right.

This is a Quillen pair because the right adjoint obviously preserves fibrations and trivial fibrations. Now we localize the projective model structure and obtain the homotopy model structure on the right-hand side. The identity functors still form an adjoint pair 
\[
\Id\colon \cal S_{\text{bifib-proj}}^{\cat A} \rightleftarrows \cal S_{\text{ho}}^{\cat A} :\! \Id,
\]
where $\sS_{\text{ho}}^{\cat A}$ denotes the homotopy model structure (which we have assumed to exist). This adjoint pair is still a Quillen pair, as a composition of the previous adjunction with the Quillen pair arising from the left Bousfield localization of the projective model structure. To show that this is a Quillen equivalence we will use \cite[Cor.~1.3.16(b)]{Hovey}. The left adjoint reflects weak equivalences between cofibrant objects, since the fibrant approximation in the homotopy model structure (approximation by the levelwise fibrant homotopy functor constructed in Proposition~\ref{homotopy-approx}) can only change the values of a bifibrant-projectively cofibrant functor in fibrant objects up to a weak equivalence.  

It remains to show that for every fibrant (homotopy) functor  $F\in \cal S_{\text{ho}}^{\cat A}$, the cofibrant replacement map $i\colon \tilde F \to F$ in the bifibrant-projective model structure $\cal S_{\text{bifib-proj}}^{\cat A}$ is  a weak equivalence in the homotopy model structure. In other words, if we apply the homotopy approximation construction $Q$ from Proposition~\ref{homotopy-approx} we obtain a projective weak equivalence. Indeed, there is a projective weak equivalence $QF\simeq F$, since $F$ is a homotopy functor. Furthermore, $Q\tilde F$ is homotopy functor bifibrant-projectively equivalent to $\tilde F$, hence also to $F\simeq QF$. So, by the 2-out-of-3 property, $Qi\colon Q\tilde F\to QF$ is a bifibrant-projective weak equivalence of homotopy functors and hence is a levelwise weak equivalence as required.
\end{proof}

\section{Motivating example} \label{motivating:example}
Before we turn to the proof of the main result, let us consider the example of ($\Delta$-generated, \cite{Fajstrup-Rosicky}) topological spaces \Top and simplicial sets $\sS$, two Quillen equivalent simplicial model categories with very different categories of small functors $\sS^\Top$ and $\sS^\sS$. For the category of functors from simplicial sets to simplicial sets we have both the bifibrant-projective model structure and the homotopy model structure constructed in the previous sections. For the case of functors from topological spaces to simplicial sets, we have several model structure to choose from. The fibrant-projective model structure is not different from the projective model structure, since all object in \Top are fibrant. The observation that simplicial functors preserve weak equivalences between cofibrant topological spaces (since every object is fibrant) suggests that for our comparison to $\sS^\sS$ we should establish the cofibrant-projective model structure on $\sS^\Top$. We do so below. Such a model structure has been previously established on the category of contravariant functors, \cite{Blanc-Chorny}, but for the category of covariant functors it is new.

\begin{proposition}
There exists the cofibrant-projective model structure on the category of small functors $\sS^\Top$, i.e., weak equivalences (resp., fibrations) are the natural transformations inducing weak equivalences (resp., fibrations) on the values of functors in cofibrant objects.
\end{proposition}

\begin{proof}
The cofibrant-projective model structure  is a particular case of the relative model structure, \cite{Chorny-relative}.  The latter exists if the solution set condition (dual to the local smallness in the case of contravariant functors) for the inclusion functor of cofibrant objects into \Top is satisfied, \cite[Prop.~2.8]{Chorny-relative}. For every uncountable regular cardinal $\kappa$ there exists a set $P_\kappa$ of $\kappa$-presentable cofibrant spaces, such that every cofibrant space is a filtered colimit of the elements of this set, \cite[Cor.~5.1]{fat}, since the domains and the codomains of the generating trivial cofibrations are finitely presentable, hence $\kappa$-presentable.

The solution set condition readily follows. Given an object $X\in \Top$, there exists an uncountable regular cardinal $\kappa$, such that $X$ is $\kappa$-presentable. Therefore, every map $X\to A$ with a cofibrant $A$ factors through the set of all possible maps $\{X\to B \,|\, B\in P_\kappa\}$.

\end{proof}

We now analyze the connection between the newly-established cofibrant-projective model structure and the homotopy model structure.

\begin{proposition}
The Quillen equivalence $(L=|-|,R=\mathrm{Sing(-)})$ between simplicial sets and topological spaces induces a Quillen equivalence $(R^*,L^*)$ between the categories of small functors $\sS^\sS$ and $\sS^\Top$ with the fibrant-projective and the cofibrant-projective model structures, respectively.
\end{proposition} 
\begin{proof}
Since $L$ preserves fibrant objects, $(R^*, L^*)$ is a Quillen pair by Proposition~\ref{Q-pair-fib-proj} between the fibrant-projective and the projective model structure, i.e., $R^*$ takes fibrations (resp., trivial fibrations) in fibrant objects into levelwise (resp., trivial) fibrations, which are also cofibrant-projective (resp., trivial) fibrations.

We will show now that $(R^*, L^*)$ is a Quillen equivalence by verifying \cite[Definition 1.3.12]{Hovey}.

Given a fibrant-projectively cofibrant $F\in \sS^\sS$ and cofibrant-projectively fibrant $G\in\sS^\Top$, consider a cofibrant-projective weak equivalence $f\colon R^*F\to G$ in $S^\Top$. The corresponding map $g\colon F\to L^*G$ in $\sS^\sS$ is constructed as a composition of the unit of the adjunction $\eta\colon F\to L^*R^*F=F(\Sing(|-|)$ with $L^*f\colon L^*R^*F\to L^*G= G(|-|)$. For every Kan complex $K\in \sS$ the natural map $K\to \Sing(|K|)$ is a simplicial homotopy equivalence (as a weak equivalence between cofibrant-fibrant objects). $F$ is a simplicial functor, hence it preserves simplicial homotopy equivalences, therefore $\eta$ is a fibrant-projective weak equivalence. The second map $L^*f$ is a projective (levelwise) weak equivalence, since $f$ is a cofibrant-projective weak equivalence and $L=|-|$ takes values in cofibrant objects. Therefore $g=L^*f\circ \eta$ is a fibrant-projective weak equivalence. 

Conversely, if we start from a fibrant-projective weak equivalence $g\colon F\to L^*G$, then the adjoint map is a composition of $R^*g$ with the counit $\epsilon\colon R^*L^*G = G(|\Sing(-)|\to G$. The first map $R^*g$ is a levelwise weak equivalence since $R=\Sing$ takes values in Kan complexes and the counit $\epsilon$ is a cofibrant-projective weak equivalence, since for every (retract of) a CW-complex $X$ the map $|\Sing(X)|\to X$ is a simplicial homotopy equivalence preserved by the simplicial functor $G$.
\end{proof}

\begin{corollary}
The homotopy model structure on $\sS^\sS$ is zig-zag Quillen equivalent to the cofibrant-projective model structure on $\sS^\Top$.
\end{corollary}

\section{Dwyer-Kan theorem for model categories}\label{DK}

In this section, we prove our main result, an extension of \cite[Thm. 2.2]{DK-diagrams} to the context of the model structures discussed above. We first prove the case where all objects are cofibrant, and then the general case. Recall that the homotopy model structure is a localization of the projective model structure.

\subsection{All objects cofibrant}
First we need to show that the adjunction $(R^*, L^*)$ is still a Quillen adjunction after the localization performed in Section~\ref{homotopy-model}.

\begin{proposition}
Consider a Quillen pair of two combinatorial model categories  $L\colon \cat A \leftrightarrows \cat B : \! R$. Then the adjunction $(R^*, L^*)$ constructed in Proposition~\ref{adjunction} between the categories of small functors equipped with the projective model structure is also a Quillen pair by Proposition~\ref{Quillen-map}. Assume in addition that all objects of \cat A and \cat B are cofibrant. Then the adjunction $(R^*, L^*)$ remains a Quillen pair for the homotopy model structure. 
\end{proposition}
\begin{proof}
By Dugger's lemma \cite[8.5.4]{Hirschhorn}, it is sufficient to verify that the right adjoint $L^*$ preserves fibrations of fibrant homotopy functors and all trivial fibrations.

Trivial fibrations are preserved since $L^*$ is a right Quillen functor in the non-localized model structure and trivial fibrations do not change (since cofibrations do not) under left Bousfield localization.

Given a fibration of two fibrant homotopy functors $f\colon F\fibr G$ in $\sS^{\cat B}$, then the induced map 
\[
L^*f\colon F(L-) = L^*F \fibr L^*G = G(L-)
\]
is again a levelwise fibration. 

Notice that $L$ preserves trivial cofibrations as a left Quillen functor. By Ken Brown's lemma, $L$ preserves weak equivalences between cofibrant objects, \cite[Corollary~7.7.2]{Hirschhorn}. Since all objects of \cat A are cofibrant, $L$ preserves weak equivalences. 

Then $L^*f$ is a fibration of homotopy functors, since $L$, $G$ and, hence, $G\circ L$ are homotopy functors, i.e., $L^*f$ is a fibration in the localized model structure.
\end{proof}

We are ready now to prove the first main result of this section stating that if the Quillen pair $(L,R)$ is a Quillen equivalence of simplicial combinatorial model categories with all objects cofibrant, then the induced Quillen pair $(R^*, L^*)$ between the categories of small functors to simplicial sets, equipped with the homotopy model structure, is also a Quillen equivalence.

\begin{theorem}\label{Quillen-equiv}
Given a Quillen equivalence $L\colon \cat A \leftrightarrows \cat B\! :R$ of two model categories with all objects cofibrant, the induced Quillen pair $(R^*, L^*)$ on the categories of small functors equipped with the homotopy model structure (obtained as a localization of the projective model structure) is also a Quillen equivalence.
\end{theorem}
\begin{proof}
We will use the criterion for a Quillen pair to be a Quillen equivalence, \cite[Corollary~ 1.3.16(c)]{Hovey}.

First we show that the right adjoint $L^*$ reflects weak equivalences of fibrant objects. Given a map of homotopy functors $f\colon F\to G$, assume that the induced map $L^*f\colon L^*F\to L^*G$ is a weak equivalence (of homotopy functors, since $L$ preserves weak equivalences).

For every  $B\in \cat B$ consider its fibrant replacement $B\trivcofib \hat B$ and put $A = R\hat B \in \cat A$. Then $LA\we \hat B$ is a weak equivalence, since $(L,R)$ is a Quillen equivalence. We obtain the following commutative diagram:
\[
\xymatrix{
F(B)
\ar[r]^{f_B}
\ar[d]_{\dir{~}}
		&	G(B)
			\ar[d]^{\dir{~}}\\
F(\hat B)
\ar[r]^{f_{\hat B}}
		&	G(\hat B)\\
F(LA)
\ar[u]^{\dir{~}}
\ar[r]
\ar@{=}[d]
		&	G(LA)
			\ar[u]_{\dir{~}}
			\ar@{=}[d]
\\
L^*F(A)
\ar[r]^{\dir{~}}
		&	L^*G(A).
}
\]

Therefore, $f_{\hat B}$ is a weak equivalence and hence $f_B$ is a weak equivalence for all $B\in \cat B$ by 2-out-of-3 property, hence $f$ is a weak equivalence.

For every cofibrant $F\in \sS^{\cat A}$, the derived unit of the adjunction $R^* \dashv L^*$ from Proposition~\ref{adjunction} is constructed as an adjoint map to the fibrant approximation in the homotopy model category $R^*F \to \cal H{R^{*}F}$:
\begin{equation}\label{to-prove}
F \to L^*\cal H{R^{*}F}.
\end{equation}
It remains to show that it is a weak equivalence in the homotopy model structure.

Note that $L^*{\cal H}(R^{*}F(-)) = L^*{\cal H} F(R(-)) = L^*\widehat{F}(R\mathrm{Fib_{\cat B}}(-)) = \widehat{F}(R\mathrm{Fib_{\cat B}}L(-))$. Since the pair $(L, R)$ is a Quillen equivalence, for all (cofibrant) $X\in \cat A$ there is a weak equivalence $X\we R\mathrm{Fib_{\cat B}}L(X)$. Hence, the initial map (\ref{to-prove}) is a weak equivalence in the homotopy model structure because we can apply $\cal H$ also to $F$ turning it into a homotopy functor. 
\end{proof}

\begin{corollary}
Assume $\cat A$ and $\cat B$ satisfy the conditions of Theorem \ref{Quillen-equiv}, and suppose that the homotopy model structures on $\sS^{\cat A}$ and $\sS ^{\cat B}$, from Theorem \ref{thm:model-for-hofunctors}, exist. Then the fibrant-projective model structures on $\sS^{\cat A}$ and $\sS^{\cat B}$ are Quillen equivalent.
\end{corollary}

\begin{proof}
By Theorem \ref{fib-proj-equiv}, the fibrant-projective model structure on $\sS^{\cat A}$ is Quillen equivalent to the homotopy model structure, and the same for $\sS^{\cat B}$. By Theorem \ref{Quillen-equiv}, the homotopy model structures are Quillen equivalent. Hence, the fibrant-projective model structures are Quillen equivalent, via a chain of Quillen equivalences (where the left adjoints are depicted): $\sS^{\cat A}_\text{fib-proj}\to  \sS^{\cat A}_\text{ho} \to \sS^{\cat B}_\text{ho}\leftarrow  \sS^{\cat B}_\text{fib-proj}$.
\end{proof}

\subsection{General case}
We adapt the definition of r-equivalences, \cite{DK-diagrams}, for simplicial model categories.
\begin{definition}\label{r-equiv}
A continuous functor $f\colon \cat A \to \cat B$ of simplicial model categories  is an r-equivalence if
\begin{enumerate}
\item
for every two bifibrant objects $A_1, A_2 \in \cat A$, the induced map $\hom(A_1, A_2)\to \hom(fA_1, fA_2)$ is a weak equivalence, and
\item
every object in the category of components $\pi_0^{\text{bifib}} \cat B$ is a retract of an object in the image of $\pi_0^{\text{bifib}}f$, i.e., for every bifibrant object $B\in \cat B$ there exists a bifibrant object $A\in \cat A$ such that $B$ is a retract of $f(A)$, up to homotopy. 
\end{enumerate}
\end{definition}

\begin{remark}
Note that  $B$ is a retract of $f(A)$, up to homotopy, if there are maps $A \overset{i}{\to} B \overset{r}{\to} A$ such that $ri\sim \Id_A$. We do not specify the kind of homotopy relation in Definition~\ref{r-equiv}(2), since for maps between bifibrant objects in a simplicial model category left, right, simplicial and strict simplicial homotopy relations coincide and are equivalence relations, \cite[9.5.24(2)]{Hirschhorn}.
\end{remark}

\begin{example}\label{main-example}
Let $\xymatrix{
\cat A
\ar@/^6pt/[r]^L_\bot		&	 \cat B
					\ar@/^6pt/[l]^R	
}$
 be a Quillen equivalence between simplicial combinatorial model categories. Put $f=\widehat L$ the composition of the left adjoint with the fibrant replacement functor in $\cat B$, and let $\widetilde{R}$ be the composition of $R$ with cofibrant replacement in $\cat A$. Then $f$ is an r-equivalence:
\begin{enumerate}
\item For all bifibrant $A_1, A_2\in \cat A$
\begin{multline*}
\hom_{\cat B}(fA_1, fA_2) = \hom_{\cat B}(\widehat L A_1, fA_2)\simeq  \hom_{\cat B}(LA_1, fA_2)\\
=\hom_{\cat A}(A_1, R\widehat LA_2)\simeq \hom_{\cat A}(A_1, A_2).
\end{multline*}
\item For all bifibrant $B\in \cat B$, factor the weak equivalence $L\widetilde R B\we B$ as a trivial cofibration followed by a fibration (also trivial by the 2-out-of-3 property): $L\widetilde R B\trivcofib \widehat L\widetilde R B\trivfibr B$. Therefore, $B$ is a retract of $\widehat L\widetilde R B$.
\[
\xymatrix{
\emptyset
\ar@{^(->}[d]
\ar[r]
			&  \widehat L\widetilde R B
				\ar@{->>}[d]^{\dir{~}}\\
B
\ar@{=}[r]
\ar@{..>}[ur]	&	B
}
\]
On the other hand, $\widehat L\widetilde R B\simeq f(\widetilde R B)$, i.e. $B$ is a retract of $f(\widetilde R B)$, up to homotopy.
\end{enumerate}
\end{example}

\begin{lemma}\label{rigid-retract}
Let $A, B\in \cat A$ be two bifibrant objects in a simplicial model category $\cat A$ such that $A$ is a retract of $B$, up to homotopy. Then there exists a bifibrant object $B'\in \cat A$ such that $B'\simeq B$ and $A$ is a strict retract of $B'$.
\end{lemma}

\begin{proof}
Suppose that the composition $A \overset{i}{\to} B \overset{r}{\to} A$ is simplicially homotopic to  the identity map on $A$: $ri\sim \Id_A$. Then $A\otimes I$ is a very good cylinder object, i.e., factors the codiagonal $A\coprod A \to A$ into a cofibration followed by a trivial fibration. Consider the following commutative diagram,
\[\xymatrix{
	& A
	  \ar[rr]^i
	  \ar@{^(->}[d]_{i_0}^{\dir{~}}	&	& B
			  		   \ar[rr]^r
					   \ar@{^(-->}[d]^{\dir{~}}	&	& A\\
A
\ar@{^(->}[r]_{i_1}^>>>>{\dir{~}}		&A\otimes I
				  \ar@/_50pt/[rrrru]_H
				  \ar@{-->}[rr]_{\imath} 	&	& B_1
				  				    \ar@{^(-->}[r]^{\dir{~}}_{i'}
				  					& B'
									   \ar@{-->>}[ur]^{r'}
}
\]
where $Hi_0=ri$ and $Hi_1=\Id_A$. Put $B_1=A\otimes I \coprod_A B$, and, in order to ensure the fibrancy of the intermediate object, we factor the natural map $B_1\to A$ as a trivial cofibration $i'$ followed by a fibration $r'$.

Hence, $A$ is a (strict) retract of the bifibrant object $B'$: $
\xymatrix{
A
\ar[r]^{i' \imath i_1}
\ar@/_10pt/[rr]_{Hi_1=\Id_A}
& B'
\ar[r]^{r'} 
& A
}
$ and $B'\simeq B$.
\end{proof}

\begin{theorem}\label{main}
Let $f\colon \cat A \to \cat B$ be an accessible functor between simplicial combinatorial model categories. Suppose $f$ preserves fibrant and cofibrant objects. Then the Quillen pair 
\[
\Lan_f \colon \sS^\cat A \leftrightarrows \sS^\cat B : \! f^*
\]
 between functor categories equipped with the fibrant-projective model structure is a Quillen equivalence if and only if $f$ is an r-equivalence of simplicial model categories.
\end{theorem}

\begin{proof}
Since $f$ preserves both fibrant and cofibrant objects, the induced adjunction $\Lan_f \dashv f^*$ is a Quillen map by Proposition~\ref{Quillen-pair-bifibrant}.

Suppose that $f$ is an r-equivalence. We will use \cite[1.3.16(c)]{Hovey} to show that  $\Lan_f \dashv f^*$ is a Quillen equivalence. In other words, we will prove that $f^*$ reflects weak equivalences of bifibrant-projectively fibrant objects and for every cofibrant $F\in \sS^\cat A$, the map 
\begin{equation}\label{QE}
F\to f^*\widehat{\Lan_f F}
\end{equation}
is a bifibrant-projective weak equivalence.

Consider a natural transformation $p\colon G\to H$ of bifibrant-projectively fibrant functors in $\sS^\cat B$. Assume that $f^*p$ is  a weak equivalence. Then for any bifibrant object $A\in \cat A$ there is a weak equivalence of simplicial sets
\[
(f^*p)(A)\colon G(f(A)) = (f^*G)(A) \we (f^*H)(A) = H(f(A))
\]

Since $f$ is an r-equivalence of model categories, Definition~\ref{r-equiv}(2) implies that for any bifibrant object $B\in \cat B$ there exists a bifibrant object $A\in \cat A$ such that $B$ is a retract of $f(A)$, up to homotopy.  By Lemma~\ref{rigid-retract} there exists $B'\in \cat B$ weakly equivalent to $f(A)$ such that $B$ is a retract of $B'$. Since $f(A)$ and $B'$ are bifibrant objects, the weak equivalence between them is a simplicial weak equivalence. Since $G$ and $H$ are simplicial functors, they preserve simplicial weak equivalences. Therefore, $p(B')$ is a weak equivalence by the 2-out-of-3 property, because $(f^*p)(A)=p(f(A))\colon G(f(A))\we H(f(A))$ is a weak equivalence by assumption. Hence, $p(B)\colon G(B)\we H(B)$ is a weak equivalence as a retract of the weak equivalence $p(B')\colon G(B')\we H(B')$. Therefore $p$ is a bifibrant-projective weak equivalence.

The retract argument implies that it is sufficient to prove (\ref{QE}) is a weak equivalence for any cellular $F\in \sS^\cat A$ with respect to the bifibrant-projective model structure. 

For every bifibrant $A\in \cat A$ 
\[
f^*\widehat{\Lan_f F}(A) = \left({\Lan_f F(f(A))}\right)^{\textup{fib}},
\]
since $f(A)\in\cat B$ is also a bifibrant object and the fibrant replacement in the bifibrant-projective model structure on $\sS^\cat B$ applies levelwise to the values of the functor in bifibrant objects.  Hence there is a bifibrant-projective weak equivalence 

\[f^*\widehat{\Lan_f F}\simeq f^*{\Lan_f F}\]

In other words, it suffices to show that the unit of adjunction 
\[
F\to f^*{\Lan_f F}
\]
is a bifibrant-projective weak equivalence.

We proceed by cellular induction. 

Suppose $F=\colim_{i < \lambda} F_i$, so that $F_0=\emptyset$ and $F_{i+1}$ is obtained from $F_{i}$ by attaching a cell 
\[
\xymatrix{
R^{A_i}\otimes \partial\Delta^n
\ar@{^(->}[d]
\ar[r]			 & F_{i}
				\ar@{^(->}[d]\\
R^{A_i}\otimes \Delta^n
\ar[r]			 & F_{i+1}\\
}
\]
if $i+1$ is a successor ordinal or $F_i=\colim_{a<i}F_a$ if $i$ is a limit ordinal. Note that $A_i\in \cat A$ is bifibrant for every $i<\lambda$. 

Since both $\Lan_f$ and $f^*$ preserve colimits, so does $f^*\Lan_f$. Since bifibrant-projective weak equivalences are preserved under sequential colimits, it suffices to show for each $i<\lambda$ that if $F_i\to f^*{\Lan_f F_i}$ is a fibrant-projective weak equivalence, then so is $F_{i+1}\to f^*{\Lan_f F_{i+1}}$. Let us consider the unit of adjunction of the homotopy pushout square above: 
\[
\xymatrix{
f^*R^{f(A_i)}\otimes \partial\Delta^n
\ar[rrr]
\ar[ddd]   &     &     & f^*\Lan_f F_i
					\ar[ddd]\\
	& R^{A_i}\otimes \partial\Delta^n
	\ar[lu]_{\dir{~}}
	\ar@{^(->}[d]
	\ar[r]			 & F_{i}
					\ar[ru]^{\dir{~}}
					\ar@{^(->}[d]\\
	& R^{A_i}\otimes \Delta^n
		\ar[ld]_{\dir{~}}
		\ar[r]			 & F_{i+1}
						\ar@{-->}[dr]\\
f^*R^{f(A_i)}\otimes \Delta^n
\ar[rrr]   &     &     & f^*\Lan_f F_{i+1}
}
\]
Then the outer square is also a pushout, moreover, this is a levelwise homotopy pushout. 

The slanted map on the right is a weak equivalence by the induction assumption.  The slanted maps on the left are bifibrant-projective weak equivalences by Definition~\ref{r-equiv}(1), since for every bifibrant $A\in \cat A$, 
\[
f^*R^{f(A_i)}(A) = \hom(f(A_i),f( A)).
\]
Hence the dashed map is also a bifibrant-projective weak equivalence as an induced map of homotopy pushouts. This completes the cellular induction proving that (\ref{QE}) is a bifibrant-projective weak equivalence, as required.

Conversely, if $f\colon \cat A \to \cat B$ preserves bifibrant objects and induces a Quillen equivalence $\Lan_f\dashv f^*$ between the categories of small functors into simplicial sets, then for any bifibrant object $A_1\in \cat A$, the induced map $R^{A_1}\to f^*\Lan_f R^{A_1} = f^* R^{f(A_1)}$ is a bifibrant equivalence, i.e., evaluating at any bifibrant object $A_2 \in \cat A$ we obtain a weak equivalence of simplicial sets
\[
\hom_{\cat A}(A_1, A_2)\, \tilde\longrightarrow \, 
\hom_{\cat B}(f(A_1), f(A_2)).
\] 
In other words, $f$ satisfies the first part of Definition~\ref{r-equiv}.  

It remains to verify \ref{r-equiv}(2), i.e., that every bifibrant object $B\in \cat B$ is a retract, up to homotopy, of $f(A)$ for some bifibrant $A\in \cat A$. 

For all (bifibrant projectively) fibrant $G\in \sS^{\cat B}$, let $\lambda$ be the maximum of the accessibility ranks of $G$ and $f$. Then
\[
G\simeq \Lan_f \widetilde{f^*G} =  \Lan_f \int^{A\in \cat A_{\lambda}} \widetilde{G(f(A))}\otimes R^A(-) = \int^{A\in \cat A_{\lambda}} \widetilde{G(f(A))}\otimes R^{f(A)}(-).
\]
Take $G= \hom(B,-)=R^B(-)$. Then
\[
\int^{A\in \cat A_{\lambda}} \widetilde{\hom(B, f(A))}\otimes \hom(f(A),-) \simeq \hom(B,-).
\]
After evaluating at $B$ and passing to connected components, we obtain a bijection, cf. \cite[Thm.~10]{weighted},
\[
\int^{A\in \cat A_{\lambda}} \pi_0\hom(B, f(A))\times \pi_0\hom(f(A),B) \cong \pi_0\hom(B,B).
\]
Let $A\in \cat A_{\lambda}$ correspond to the identity on the right hand side. Then $B$ is a retract, up to homotopy, of $f(A)$.
\end{proof}

\section{Invariance of Goodwillie calculus under Quillen equivalence}\label{Goodwillie}
Given a Quillen equivalence $f$ of simplicial combinatorial model categories, consider the model categories of homotopy functors (bifibrant-projective model structure on the categories of small functors) from this Quillen pair to simplicial sets. This is a starting point for Goodwillie's calculus of homotopy functors, \cite{Goo:calc3, Kuhn:overview}. Consider the localization of these categories such that the fibrant objects are the fibrant $n$-excisive functors. Does $f$ induce a Quillen equivalence of these model structures? In other words, is Goodwillie calculus invariant under Quillen equivalence? Under a few additional conditions the answer is `yes.'

\begin{theorem} \label{thm:goodwillie}
Let $\xymatrix{
\cat A
\ar@/^6pt/[r]^L_\bot		&	 \cat B
					\ar@/^6pt/[l]^R	
}$
 be a Quillen equivalence between simplicial combinatorial model categories.  Put $f = \widehat L$. Then the Quillen pair of functor categories 
\[
\Lan_f \colon \sS^\cat A \leftrightarrows \sS^\cat B : \! f^*
\]
is a Quillen equivalence. Moreover, if we left Bousfield localize the functor categories so that the local objects are the $n$-excisive functors in the bifibrant-projective model structures on both sides, then the adjunction $\Lan_f  \dashv f^*$ is a Quillen equivalence of the model categories of $n$-excisive functors.
\end{theorem}

\begin{proof}
That the first Quillen pair is a Quillen equivalence follows from Theorem \ref{main} and Example \ref{main-example}. In order to conclude that the categories of $n$-excisive functors are Quillen equivalent we will apply Theorem~\ref{Q-equiv-after-loc}.

First, $f^*$ preserves $n$-excisive functors. This is because $f = \widehat L$ preserves homotopy pushouts of cofibrant objects, hence $f$ preserves strongly cocartesian cubes.

It remains to show that $\Lan_f$ commutes with the $n$-th polynomial approximation functor. Let us denote by  $P_n^{\cat A}$ and $P_n^{\cat B}$ the $n$-excisive approximations of functors with domain in \cat A and \cat B respectively. Then we need to prove that for all cofibrant $F\in \sS^{\cat A}$ there is a bifibrant-projective weak equivalence of functors: $\Lan_f P_n^{\cat A} F \simeq P_n^{\cat B}\Lan_f F$.

We next prove that $\Lan_f$ takes $n$-excisive functors to $n$-excisive functors. We will need the right adjoint $R$ to preserve homotopy pushout squares of bifibrant objects to verify this property. This is true, in turn, since $R$ is a part of a Quillen equivalence, hence its total derived functor is an equivalence of homotopy categories. It follows that there is an equivalence of homotopy categories of diagrams indexed by the category $\bullet \leftarrow \bullet \rightarrow \bullet$. Since the homotopy pushout is a left adjoint to the constant functor, it is preserved by any equivalence of categories, so $R$ preserves homotopy pushout squares.

Let $F$ be a cofibrant functor. Suppose that $P_n^{\cat A}F$ is $\lambda$-accessible functor. Let $\cat A_{\lambda}\subset \cat A$ be the subcategory of $\lambda$-presentable objects, then
 \begin{align*}
& \Lan_f P_n^{\cat A}F = \Lan_f \int^{A\in \cat A_{\lambda}}\hom(A,-)\otimes P_n^{\cat A}F(A) \\
&= \int^{A\in \cat A_{\lambda}}\Lan_f \hom(A,-)\otimes P_n^{\cat A}F(A) = \int^{A\in \cat A_{\lambda}} \hom(f(A),-)\otimes P_n^{\cat A}F(A)\\
 & = \int^{A\in \cat A_{\lambda}} \hom(\widehat L(A),-)\otimes P_n^{\cat A}F(A) \simeq \int^{A\in \cat A_{\lambda}} \hom(L(A),-)\otimes P_n^{\cat A}F(A)\\
 & =\int^{A\in \cat A_{\lambda}} \hom(A,R(-))\otimes P_n^{\cat A}F(A) = R^*P_n^{\cat A}F.
 \end{align*}
Since $R$ preserves homotopy pushouts, it preserves strongly cocartesian cubes, hence $R^*P_n^{\cat A}F \in \sS^{\cat B}$ is an $n$-excisive functor.

Consider the diagram
\[
\xymatrix{
\Lan_f F 
\ar[d]
\ar[r]		& P_n^{\cat B} \Lan_f F 
			\ar@{..>}@/^/[dl]	\\
\Lan_f P_n^{\cat A} F 
\ar@{-->}@/^/[ur]	
}
\]
The dotted arrow exists by the initial, up to homotopy, property of $n$-excisive approximation in $\sS^{\cat B}$, since $\Lan_f P_n^{\cat A} F$ is $n$-excisive.

The dashed arrow exists by the same reason in the adjoint diagram:
\[
\xymatrix{
F 
\ar[d]
\ar[r]		& f^*P_n^{\cat B} \Lan_f F 
				\\
P_n^{\cat A} F.
\ar@{-->}[ur]	
}
\]
Both maps are unique, up to homotopy, hence mutually homotopy inverse. Therefore, $\Lan_f$ commutes with polynomial approximation, up to homotopy. 

We have verified the conditions of Theorem~\ref{Q-equiv-after-loc}, so the Quillen pair $\Lan_f\vdash f^*$ is a Quillen equivalence between the categories of small functors equipped with the $n$-excisive model structures.
\end{proof}

\begin{example}
For the Quillen equivalence $|-|\colon \sS \leftrightarrows \Top :\! \Sing$  the categories of $n$-excisive functors $\sS^{\sS}$ and $\sS^{\Top}$ are Quillen equivalent. In other words, simplicial sets and topological spaces have the same calculus of homotopy functors.
\end{example}

\appendix

\newcommand{\FA}{\ensuremath{\mathcal{F}_{\cat A}}}
\newcommand{\FB}{\ensuremath{\mathcal{F}_{\cat B}}}

\section{Localization of a Quillen equivalence}\label{compare-localizations}

Given two Quillen equivalent model categories, consider a left Bousfield localizations of both sides. Under what conditions are the resulting localized categories Quillen equivalent again? We provide a new answer, needed for our Theorem \ref{thm:goodwillie}, in the following theorem. We fix the following notation: for a model category \cat A equipped with the homotopy localization functor $\FA^{-1}\colon \cat A\to \cat A$ the left Bousfield localization of \cat A with respect to  $\FA^{-1}$-equivalences is denoted by $\FA^{-1}\cat A$. For all $A\in \cat A$ we write $\FA^{-1}(A)$ for the fibrant replacement of $A$ in $\FA^{-1}\cat A$ in order to distinguish it from the fibrant replacement in \cat A.

\begin{theorem}\label{Q-equiv-after-loc}
Let $\xymatrix{
\cat A
\ar@/^6pt/[r]^L_\bot		&	 \cat B
					\ar@/^6pt/[l]^R	
}$
 be a Quillen equivalence between simplicial model categories. Suppose there exist left Bousfield localizations $\FA^{-1}\cat A$ of $\cat A$ and $\FB^{-1}\cat B$ of $\cat B$ such that 
 \begin{enumerate}
 \item
 $R$ takes $\cal F_{\cat B}$-local objects to  $\cal F_{\cat A}$-local objects;
 \item
 $L$ commutes with the localization, i.e., for all cofibrant $A\in \cat A$ the map $L\FA^{-1} (A) \to \FB^{-1}(LA)$ is a weak equivalence. The latter map is adjoint to the lift in the following commutative square in $\FA^{-1}\cat A$:
 \[
 \xymatrix{
 A
 \ar[r]
 \ar@{^(->}[d]^{\dir{~}} & RLA \ar[r] & R\FB^{-1}(LA)
 							\ar@{->>}[d]\\
 \FA^{-1}(A) 
 \ar[rr] 
 \ar@{-->}[urr]& & \ast
 }.
 \]
 \end{enumerate}
 Then there exists a Quillen equivalence of the localized model categories 
 $$
 \xymatrix{
\FA^{-1}\cat A
\ar@/^6pt/[r]^L_\bot		&	 \FB^{-1}\cat B
					\ar@/^6pt/[l]^R.
}
$$

\end{theorem}
\begin{proof}
$L$ is a left Quillen functor between \cat  A and \cat B, hence it preserves the cofibrations in the localized model structures as well. In order to show that $L$ remains a left Quillen functor after the localization we need to verify that for every cofibration $A_1\cofib A_2$, which is also an $\FA$-local equivalence, the cofibration $LA_1\cofib LA_2$ is an $\FB$-local equivalence. In other words, we need to prove that $\FB^{-1}LA_1\to \FB^{-1}LA_2$ is a weak equivalence. Since $L$ commutes with the localization, it suffices to show that $L\FA^{-1}A_1\to L\FA^{-1}A_2$ is a weak equivalence, which readily follows from the assumption that $A_1\cofib A_2$ is an $\FA$-local equivalence and the fact that $L$ preserves weak equivalences of cofibrant objects.

We will use  \cite[1.3.16(c)]{Hovey} to show that $L \dashv R$ is a Quillen equivalence.
\begin{enumerate}
\item
$R$ reflects local equivalences of local objects, since these are just weak equivalences in $\cat A$ and $\cat B$, and $R$ reflects weak equivalences.
\item
For every cofibrant $A\in \cat A$ we need to show that the map $A\to R\FB^{-1} LA$ is an $\FA^{-1}$-local equivalence.  We need to rely on the assumption that $L$ commutes with the localization, i.e., that $L\FA^{-1} A \to \FB^{-1}LA$ is a weak equivalence in \cat B.

Consider the commutative square in the model category $\FB^{-1}\cat B$:
\[
\xymatrix{
{L\FA^{-1} A}
\ar@{^(->}[d]^{\dir{~}}
\ar[r]^{\dir{~}} 				& {\FB^{-1}LA}
							\ar@{->>}[d]\\
{\widehat{L}\FA^{-1} A}
\ar[r]
\ar@{-->}[ur]_h						& \ast
}
\]
The lift $h$ is a weak equivalence of fibrant objects in $\FB^{-1}\cat B$. Hence $Rh$ below is a weak equivalence:
\[
\xymatrix{
A
\ar@{^(->}[d]_{\txt{\tiny loc. equiv.}}
\ar[r]						& R{\FB^{-1}LA}
							\\
\FA^{-1} A
\ar[r]^>>>>>{\dir{~}}
						& R{\widehat L\FA^{-1} A}
							\ar[u]_{Rh}^{\dir{~}}
}
\]
In the diagram above the lower horizontal map is a weak equivalence since $L\dashv R$ is a Quillen equivalence between \cat A and \cat B. Therefore, the map $A\to R\FB^{-1} LA$ is a local equivalence in $\FB^{-1} \cat A$.
\end{enumerate}
\end{proof}

\begin{example}
Consider the Quillen equivalence $|-|\colon \sS \leftrightarrows \Top :\! \Sing$.  And consider the localization of both sides with respect to integral homology. Then the conditions of Theorem~\ref{Q-equiv-after-loc}  are readily verified, hence the category of $H\ZZ$-local topological spaces is Quillen equivalent to the category of $H\ZZ$-local simplicial sets.
\end{example}

\bibliographystyle{abbrv}
\bibliography{../Xbib}


\end{document}